# Analyse du travail des étudiants de l'ENEPS sur la plateforme de cours en ligne

Gilles Aldon, Corinne Raffin, équipe EducTice, IFÉ-ENS de Lyon

Résumé

L'objet de ce rapport est de décrire le travail réalisé pour la conception des modules d'un SPOC (Small Private On line Course) et l'analyse de leurs effets sur les apprentissages des étudiants. Nous avons proposé une méthodologie d'évaluation de cette plateforme en alliant une approche didactique dans la conception et l'analyse des contenus et une approche ergonomique.

Abstract

The purpose of this report is to describe the work done to design the modules of a SPOC (Small Private On line Course) and analyse their effects on students' learning. We proposed a methodology for evaluating this platform by combining a didactic approach in the design and analysis of content and an ergonomic approach.

# Introduction

Depuis trois ans nous avons participé à l'analyse et à la réalisation des modules de la plateforme de mise à niveau des étudiants de l'ENEPS en mathématiques. Le point de départ de ce travail repose sur le constat des difficultés rencontrées par certains étudiants à suivre les cours de première année du fait de concepts mathématiques ou de techniques mal ou pas maîtrisés. L'objectif de cette plateforme est ainsi de permettre à ces étudiants de tester leurs connaissances mathématiques et de leurs permettre de comprendre ces notions mathématiques suffisamment pour qu'elles ne soient pas un obstacle à la compréhension des cours de physique.
Les modules ont été pensés sous un même format :

- Présentation du module : une première accroche problématique pour faire percevoir aux étudiants la pertinence du concept ou de l'outil qui est travaillé dans le module.

- Des cours en vidéo permettant de détailler les contenus qu'il s'agit de maîtriser.

- Des exercices pour s'entraîner corrigés

- Une évaluation "intermédiaire" sous forme de quiz permettant aux étudiants d'auto-évaluer leurs connaissances

- Des exercices supplémentaires pour s'entraîner

- Une évaluation finale sous forme de quiz.

Après que la majorité des modules aient été implémentés, nous avons voulu observer le comportement des étudiants confrontés à la plateforme et recueillir leurs avis quant à l'utilisation qu'ils en ont fait.
Nous décrivons dans la suite la méthode utilisée reposant tout d'abord sur une analyse par inspection de la plateforme puis sur une analyse empirique avec les observations et entretiens réalisés avec les étudiants. Les cadres d'analyse utilisés sont d'une part l'analyse issue de l'ergonomie avec les concepts d'utilité, d'utilisabilité et d'acceptabilité (Tricot, 2003) et de la théorie des situations didactiques (Brousseau, 1986) mais aussi de la théorie anthropologique du didactique (Chevallard, 1985).
Nous donnons dans une première partie les résultats essentiels de l'analyse sous forme de 6 remarques que nous développons et expliquons dans la seconde partie.

# Résultats des observations en 6 remarques

Ces remarques s'appuient sur les analyses de la plateforme en tant qu'outil (analyse par inspection) et les observations et entretiens réalisés avec les étudiants (analyse empirique). Les termes utilisés seront explicités dans la partie "méthodologie" et ces remarques seront développées dans la partie d'analyse.

## Remarque 1

La plateforme dans sa diversité d'approches (cours, exercices, auto-évaluation, explication plus particulière, évaluation) est suffisamment **flexible** pour s'adapter aux besoins très divers des étudiants.
Du point de vue des contenus des modules :
- la possibilité de voir des vidéos avec des explications données, ou de télécharger les cours, ou les corrections d'exercices,
- la possibilité de s'auto-évaluer avec des QCM intermédiaires

participent à la flexibilité de la plateforme
Du point de vue des besoins des étudiants
- les exercices nombreux,
- les QCM adaptés

permettent aux étudiants d'adapter leur apprentissage en fonction de leurs connaissances initiales.

## Remarque 2

Les modules sont **utiles** pour les étudiants pour :
1. se conforter dans les connaissances nécessaires pour la suite de leur intégration en leur donnant un point de référence initial.
2. Pour revoir quelques notions qui sont fondamentales dans la suite de leurs cursus.

## Remarque 3

Il apparaît que les vidéos sont **utiles** pour les étudiants mais que les utilisations qu'ils en font sont diverses et adaptées à leurs stratégies d'apprentissage :
   - lecture continue de la vidéo, ce qui souvent provoque le sentiment de longueur

- lecture avec des arrêts pour vérifier ou faire les exercices
  - lecture avec des accélérés au moment des parties supposées connues

Il est donc intéressant de garder ce format de vidéos qui permet à chacun d'utiliser ce support de cours à sa façon.

# Remarque 4

Les exercices sont adaptés et fondamentaux pour l'entraînement des étudiants. Même lorsque le sujet paraît simple, c'est à travers la réalisation effective de l'exercice qu'ils confortent leur assurance à savoir résoudre des problèmes. D'une façon générale, les exercices, les QCM et les exemples utilisés dans les modules sont **acceptables** par les étudiants qui les ont expérimentés.

# Remarque 5

L'ergonomie générale rend la plateforme **utilisable** et n'est pas un frein à l'acceptabilité des contenus. Cependant, les observations et les entretiens montrent que des progrès pourraient être faits pour que l'accès aux modules soit simplifié.

# Remarque 6

Les étudiants sous-estiment l'importance des connaissances contenues dans les modules. En particulier, ils sous-utilisent la capitalisation des travaux réalisés dans ces modules pour en faire une référence initiale pour les disciplines qu'ils vont rencontrer. Il pourrait être intéressant de demander aux étudiants de constituer un port-folio pour mettre en évidence les connaissances acquises mais aussi les difficultés rencontrées.

# Méthodologie

Dans un premier temps nous considérerons les dimensions d'analyse d'un EIAH d'une part pour justifier les choix des contenus et de la structure et d'autre part pour analyser l'utilisation de la plateforme par les étudiants de l'ENEPS. Les cadres théoriques de didactique des mathématiques s'appuient sur la théorie des situations didactiques (TSD) (Brousseau, 1986) et sur la théorie anthropologique du didactique (TAD) (Chevallard, 1984). La plateforme propose aux étudiants des modules d'apprentissage pendant lesquels ils sont confrontés à un environnement particulier qui, de part ses rétroactions, permet d'apprendre les contenus mathématiques proposés par les auteurs. Les questions des apprentissages effectifs sont bien sûr posées qui mettent en avant à la fois les contenus et la pertinence de ces contenus pour l'apprentissage des mathématiques mais aussi les institutionnalisations qui pourront en être faites : qu'est ce que les étudiants pourront considérer comme acquis une fois les modules conduits ? Ces questions sont examinées en prenant comme base les observations et les entretiens avec les étudiants. Les exercices proposés, tout comme les activités des différents modules ont été analysés en s'appuyant sur une analyse praxéologique des contenus. Pour Chevallard (ibid.) toute activité humaine confronté à un certain type de tâches peut être décrite en distinguant une partie technique de réalisation de la tâche et une partie conceptuelle d'explicitation de la technique utilisée. Une praxéologie est ainsi décrite, à l'intérieur d'une institution donnée, comme un quadruplet (t, T, tau, Theta) où t est un type de tâche, T une technique pour réaliser la tâche particulière envisagée, tau une technologie, c'est à dire un discours justifiant cette technique et Theta la théorie sous-jacente à cette justification. Pour les contenus mathématiques des modules, la technique de réalisation des tâches proposées est importante mais nous avons toujours essayé de proposer des justifications des techniques, en particulier en montrant l'importance du type de tâche proposé dans un contexte plus large. L'exemple développé ci-dessous dans l'encadré 2 montre comment les différents contenus peuvent être construits et pourront l'être dans une perspective d'élargissement des modules.
En croisant les apports que ces cadres théoriques peuvent amener à la réflexion et à l'analyse, nous analysons la plateforme comme un EIAH spécifique dont les contenus mathématiques constituent l'objectif d'apprentissage.

La méthodologie d'observation s'est déroulée en trois phases distinctes :

1. une analyse des log des étudiants pour constituer des types de comportements sur la plateforme,

2. une observation d'étudiants choisis sur cette typologie en condition de travail,

3. un entretien sous forme de q-sorting.

# Analyse par inspection

Une analyse par inspection consiste à parcourir la plateforme en appliquant une grille d'observation provenant du cadre théorique considéré. En reprenant les formulations de Tricot (2003), nous évaluons la plateforme selon les trois dimensions de l'utilité, l'utilisabilité et l'acceptabilité.

## Utilité

La construction de la structure des modules s'est appuyée sur la volonté de mettre en adéquation d'une part l'objectif défini et l'apprentissage effectif et d'autre part l'adéquation entre le dispositif et les connaissances à acquérir. Ainsi l'objectif de maîtrise de techniques calculatoires (développement, factorisation, résolution d'équations, parenthésages, etc.) résulte de la structure des modules conçus dans la perspective d'un entraînement suffisant aux techniques à maîtriser. Mais aussi en donnant une justification de ces techniques par un exemple introductif suffisamment probant pour que les étudiants comprennent le sens des calculs à apprendre. Ainsi par exemple, parler d'exponentiation passe à la fois par la maîtrise des règles usuelles de calcul mais aussi par la compréhension de cette opération et son lien avec les opérations usuelles ; l'exemple de la feuille pliée (encadré 1.) est paradigmatique de cette compréhension des objets sur lesquels les étudiants vont ensuite entraîner leur compétences calculatoires. Le résultat paradoxal montre bien la croissance très rapide des puissances d'un nombre et donne une idée de la différence de l'exponentiation et de la multiplication ou de l'addition.

---

Une feuille de papier de 80g a une épaisseur de 0,1mm. On la plie en deux, puis on plie la feuille pliée encore en deux, etc. Combien de fois devriez vous la plier en deux pour atteindre la hauteur de la Tour Eiffel (312m=312 000 mm) ?
Répondez à la question de façon intuitive, sans faire de calculs
a) 22 fois b) 100 fois, c) 3120000 fois d) beucoup plus !

---

Encadré 1 : Introduction au module Puissances

Par ailleurs, les tâches proposées prennent en compte les dimensions techniques et technologiques des praxéologies liées aux types de tâches du module. Ainsi, les différentes compétences nécessaires à la compréhension et la maîtrise des contenus ont été analysées et participent à la diversification des tâches. Par exemple, l'analyse des compétences travaillées dans ce même module conduit à la proposition des différents exercices et en particulier la remarque finale a permis l'élaboration de la situation d'introduction.

> Si on essaye de considérer les types de tâches qui peuvent apparaître dans l'étude des puissances on peut les décrire comme :
>
> T0 : Transformer un nombre en écriture sous forme d'une puissance en une écriture sous forme décimale.
> T'0 : Déterminer le signe d'une expression comportant des puissances.
> T1 : Transformer un nombre décimal bien choisi sous la forme d'une écriture avec une puissance.
> Exemple : 125
> T2 : Écrire un nombre décimal non nul en écriture scientifique.
> T3 : Écrire un nombre en écriture scientifique sous une forme décimale.
> T4 : Écrire un nombre écrit avec des puissances de 10 en notation scientifique.
> **Pour le calcul :**
> T5 : Calculer un produit de puissances.
> T6 : Calculer un quotient de puissances.
> T7 : Calculer une puissance de puissances.
> **Pour les problèmes**
> T8: Interpréter et reconnaître une situation qui relève d'un modèle multiplicatif répété par une puissance\footnote (Guerrin, 2012)
> Le type de tâche T'0 est ajouté à la publication.
> Tout ce qui relève de l'écriture scientifique ne concerne pas ce cours (T2, T3, T4).
> Exercice 1 T1
> Exercice 2 T0
> Exercice 3 T5, T6
> Exercice 4 T1, T7
> Exercice 5 T'0
> Exercice 6 T5, T6, T7
> T7 n'apparaît pas beaucoup (1 fois dans l'exercice 3 et 2 fois dans l'exercice 6 et pourtant c'est sûrement une grande difficulté dans les calculs.)
> Le seul type de tâche non représenté est T8 qui donne du sens à la notion de puissance et qui pourrait permettre de mieux intégrer les propriétés des puissances. Il me semble que ce serait bien de penser un exercice issu du domaine de compétence des étudiants mettant en jeu un principe d'itération de multiplications.

Encadré 2 : un exemple d'analyse pour la construction d'un module

## Utilisabilité

Les contraintes de l'université de Grenoble font que la plateforme utilisée a été la plateforme Chamilo développée et utilisée par l'université Joseph Fourier. Par conséquent, la création des cours a dû être faite en tenant compte de ces contraintes et sans pouvoir agir sur la gestion des erreurs, sur la facilité de navigation et sur l'ergonomie même des cours proposés.

## Acceptabilité

La forme même donnée aux cours et les incitations envoyées aux étudiants pour travailler sur la plateforme est *a priori* acceptable pour les étudiants. Dans les logs de la plateforme, on peut voir les temps de connection, les heures d'entrée et de sortie et les parties des modules qui ont été visitées. Ces éléments permettent de construire des profils d'utilisation de la plateforme et constitue une première approche de l'étude de l'acceptabilité. Les connections nombreuses, les quiz réalisés, la motivation des étudiants tendent à montrer que, du point de vue des utilisateurs, la plateforme est acceptable.
La volonté de l'institution de créer, de suivre, d'alimenter et d'inciter les étudiants à l'utiliser est du point de vue de l'institution un gage d'acceptabilité. Ces premiers éléments doivent être confirmés et vérifiés dans la partie suivante.
Le choix des observations a été réalisé sur les observations des log des étudiants et des travaux dont les logs peuvent témoigner. Les étudiants "sélectionnés" pour les observations et les entretiens l'ont été sur la base des remarques faites dans le tableau suivant. Des étudiants sont dans les deux sections concernés par cette étude, en rouge, les étudiants en génie mécanique, en bleu, les étudiants en génie civil.

| étudiants | durée | nombre modules -réussites | organisation de W |
|---|---|---|---|
| BM | 0h56 | module 1<br>50 à 100% | 1 tentative<br>très peu de temps sur vidéo<br>mais du temps sur pdf cours |
| OC | 5H40 | 3iers | module 1 : laborieusement avec bcq de tentatives aux tests<br>module 2 : très rapidement |
| MC | 11H50 | tout<br>score > 75 | n'a pas cherché à améliorer ses scores lorsque < 50 |
| DV | 1H23 | module 1<br>score de 50 à 100 | une tentative<br>peu de temps sur vidéo<br>plus sur pdf cours |
| LE | 3H03 | 4 modules<br>scores<75 | parfois a cherché à améliorer ses scores d'autres fois non<br>de 9mn à 1H36 sur un module |
| FL | 1H38 | 3 modules | peu ou pas de temps sur vidéos |

|  |  | scores<75 | plusieurs tentatives aux tests<br>ne fait toujours le QCM final |
| --- | --- | --- | --- |
| NM | 4H08 | 3 modules<br>scores <75 | plusieurs tentatives aux tests |
| TE | 5H | tout | n'a pas toujours fini<br>une seule tentative aux tests même si <75<br>+ ou - sérieusement selon les modules |
| LL | 8H49 | tout | scores irréguliers<br>plusieurs tentatives(2) parfois<br>semble passer du temps sur les exos |
| BY | 2H38 | 4 modules | une tentative même si <75<br>zappe les vidéos |
| BJ | 2H45 | tout<br>score >75 sauf module 3 | semble consacrer du temps à ce dont il a besoin<br>directement aux QCM<br>n'améliore pas toujours ses scores |
| GL | 6H44 | tout<br>scores >62 | plusieurs tentatives jusqu'à scores > 62,5<br>n'a pas toujours fait le QCM final |

## Flexibilité

"*The criterion flexibility refers to the means available to the users to customize the interface in order to take into account their working strategies and/or their habits and the task requirement. Flexibility is reflected in the number of possible ways of achieving a given goal*" (Bastien & Scapin, 1993, page 30).

La flexibilité peut être considérée du point de vue de l'EIAH lui même et en ce sens renvoie à l'utilisabilité ou bien du point de vue des possibilités offertes aux utilisateurs pour adapter l'environnement à leurs propres intentions d'apprentissage. Il est en particulier intéressant de voir en quoi la plateforme permet aux étudiants d'accomplir une tâche de différentes façons en maintenant une cohérence avec les intentions des auteurs. Et de la même façon, en quoi les modules proposés permettent à tous les étudiants de trouver des éléments suffisants pour s'adapter à leurs apprentissages.

## En conclusion

L'utilité de la plateforme est attestée à la fois par l'analyse des enseignants de l'ENEPS et la mise en adéquation des problèmes rencontrés dans les enseignements de première année. La conception de la structure de la formation proposée sur la plateforme permet d'insister sur cette adéquation et confirme l'utilité de la plateforme dans cette première approche.

L'utilisabilité est plus difficile dans cette partie à analyser du fait des contraintes imposées par la plateforme utilisée. L'analyse empirique permettra de donner des éléments probants pour cette propriété importante.

Enfin, l'acceptabilité, vue à la fois du côté des étudiants que du côté des enseignants est attesté par les premières utilisations et sera complétée par l'étude empirique.

## Analyse empirique

Cette analyse se fonde sur des observations en situation ainsi que des entretiens conçus comme q-sorting interviews.

Observations : nous avons demandé à un groupe d'étudiants de travailler sur un module qui n'avait pas été encore mis en place sur la plateforme. Deux observateurs ont observé les comportements et le travail des étudiants. Des écrans d'ordinateurs ont été enregistré ainsi que les dialogues des étudiants.

Q-sorting interview : il s'agit de proposer aux étudiants des avis sous forme de phrases exprimant une position tranchée sur les sujets à débattre et de demander à un groupe de se mettre d'accord sur le degré d'acceptation de cette phrase : Tout à fait d'accord, d'accord, plutôt pas d'accord, pas du tout d'accord. Cette phase de travail de groupe a été enregistrée (son et images).

### Observations

La grille suivante permet de mettre en évidence les comportements des étudiants en situation par rapport à leur profil étudié précédemment. Pour chaque étudiant "sélectionné" une observation plus fine de l'utilisation des vidéos du temps passé sur chaque section du cours permet d'affiner ce profil.

| Etudiants | Observations faites sur les modules réalisés avant l'observation | | | Observations sur place | | | |
|---|---|---|---|---|---|---|---|
| | durée | nombre modules-réussites | organisation de W | Démarche | Outils numériques utilisés | Aide extérieure à la plateforme | Vidéos-scores |
| BM | 0h56 | module 1 50 à 100% | 1 tentative très peu de temps sur vidéo mais du temps sur pdf cours | Ds l'ordre<br>- présentation : 7mn<br>- Plein écran<br>- Fait les calculs ds sa tête<br>- QCM final | La souris pour compter le déplacement de la virgule | Entre aide avec valentin Domois | Tps vidéos : 2 à 3mn30<br>QCM inter :5/5<br>QCM final : 1/4 |
| DV | 1H23 | module 1 score de 50 à 100 | une tentative peu de temps sur vidéo plus sur pdf cours | Ds l'ordre<br>- présentation : non<br>- Regarde attentivement et rapidement les vidéos<br>- Ne clique jamais<br>- Exos Pdf avec papier de calepin<br>- QCM inter mais choqué par ses résultats | Pdf | Entre aide avec Bourgin mathias<br><br>Petit papier | Tps vidéos :0 à 5 mn<br>QCM inter : 2/5<br>QCM final : 1/4 |
| NM | 4H08 | 3 modules scores <75 | plusieurs tentatives aux tests | Ds l'ordre<br>- présentation : 4mn31<br>- QCM inter avec calculatrice et/ou ds sa tête | | Calculatrice | Tps vidéos : 1 à 5 mn<br>QCM inter :2/5 5/5<br>QCM final:1/4 |
| PD | 0h23 | 2modules Scores <75 | Ne fait que les tests | - plein écran<br>- ds l'ordre<br>- présentation : 6mn37<br>-vidéo rapidement<br>- ne clique jamais<br>- QCM inter : 2 essais / déplace la souris pour lire / lit résultats et réponses<br>- lit pdf exos et cours<br>- QCM final : 4 essais | Pdf<br><br>La souris pour lire | | Tps vidéos : 2-3 mn<br>QCM inter :2/5 4/5<br>QCM final : 1/4 1/4 2/4 1/4 |
| KA (vidéo) | 0h52 | 5 modules Score de 25 à75 | 1 tentative pour tests intermédiaires 2 pour test final Jette un œil rapide sur tout | - plein écran<br>-Ds l'ordre<br>- présentation :9 mn<br>- prend des notes sur mise en train<br>- utilise les arrêts clics et fait calculs par écrit ou les recopie<br>- QCM inter : utilise sa calculatrice. 1essai<br>- cherche les exos Pdf sur papier<br>- QCM final 1/4 | Arrêts clics<br><br>Pdf | Support papier<br><br>Entre aide avec Emeline (millionième) | Tps vidéos : 3-4mn<br>QCM inter : 5/5<br>QCM final : 1/4 |
| LE | 1h10 | 5 modules Score de 37 à 100 % | 1 ou 2 tentatives pour QCM inter. Jusqu'à 4 tentatives pour QCM final | - plein écran<br>- Ds l'ordre<br>- présentation : 13s<br>- papier devant elle : prend des notes, recopie les règles de calcul<br>- QCm inter. | Arrêts sur image avec barre d'espace<br><br>Google : recherche sur millionième | Support papier<br><br>Calculatrice<br><br>Entre aide avec Ali recherche sur millionième | Tps vidéos : 2 à 7 mn<br>QCM inter : 3/5 3/5<br>QCM final : 1/4 1/4 2/4 2/4 |
| FL | 1H38 | 3 modules scores<75 | peu ou pas de temps sur vidéos plusieurs tentatives aux tests ne fait tjrs le QCM final module 6 : intro 6s | - commence par QCM intermédiaire<br>- papier sur la table + calculatrice<br>- QCM final fait | calculatrice | Calculatrice (vérification calculs) | Tps vidéos : 2 à 1 mn20<br>QCM inter : 4/5<br>QCM final : 2/4 2/4 3/4 2/4 4/4 |
| LE | 3H03 | 4 modules scores<75 | parfois a cherché à améliorer ses scores d'autres fois non de 9mn à 1H36 sur un module module 6 : intro 9mn33s | - visionne les vidéos cours en premiers<br>- semble faire le calcul en suivant (mimant l'écriture?) avec la souris<br>- calculs à la main sur papier | Vidéo casque audio souris pour suivre le calcul. | | Tps vidéos : 2 à 4 mn<br>QCM inter : 1/5 3/5<br>QCM final : Non |
| MA | 3h14 | 5 modules Scores de 66 à 100 | Plusieurs tentatives lorsque besoin Très rapide sur les contenus module 6 : intro 2mn49s | - vidéo en plein écran<br>- QCM intermédiaire q1 : faux reste un moment sur la question<br>- QCM : 6,4 $10^7$ → 640 000 000 reste un moment pour comprendre son erreur.<br>- revient au cours<br>- QCM final (un pourcentage de 93 % apparaît alors que 3 questions sur 4 sont réussies) | Vidéo casque audio | Entraide avec NJ | Tps vidéos : 7s à 5 mn<br>QCM inter : 4/5 5/5<br><br>QCM final : 2/4 4/4 |
| NJ | 4h19 | 5 modules Scores >75 | 1 tentative module 6 : intro 5mn54s | - regarde exercice corrigé<br>- regarde le cours vidéo<br>- utilise sa calculatrice pour faire le QCM intermédiaire (il teste toutes les réponses) | Pdf vidéos casque audio TI Nspire CAS | Entraide avec MA | Tps vidéos : 2s à 7 mn<br>QCM inter : 5/5<br>QCM final : 1/4 1/4 |
| NT | 0h36 | 4 modules Scores de 44 à 66 | 1 à 3 tentatives Tout à une vitesse grand V module 6 : intro 9s | - visionne les vidéos cours en premier<br>- sur la page exercice (sans papier)<br>- QCM final : 2 exacts, 2 faux : refait les calculs avec la calculatrice → refait le pdf sur l'écran calculatrice de l'ordinateur | QCM Vidéos casque audio | Calculatrice de l'ordinateur | Tps vidéos : 1s à 2 mn46<br>QCM inter : 3/5<br>QCM final : 3/4 |
| OC | 5H40 | 3iers Score de 42 à 100 | module 1 : laborieusement avec bcq de tentatives aux tests<br>module 2 : très rapidement<br>module 6 : intro 3s | Se connecte et trouve la bonne plateforme laborieusement ⇒ 6mn 49<br>Zappe la présentation et va directement au 1ier cours<br>Notations : branche son casque / entièrement sans interruption<br>Multiplications et divisions par 10 : change la luminosité de son écran / entièrement sans interruption<br>Rappels, calculs avec les puissances : accélère les exemples et ne regarde pas tout<br>Ecriture scientifique : la moitié, pas les exemples<br>Comparaison : accélère et ne regarde pas tout<br>Exercice corrigé 1 : ne regarde pas tout<br>Exercice corrigé 2 : ne regarde qu'une partie du premier exemple<br>Exercice corrigé 3 : en partie en accélérant | Pdf | Aide de camarades QCM final | Tps vidéos : 30s à 4 mn 30<br>QCM inter :1/5 3/5<br>QCM final : q1 : 1mn17 avec aide du voisin / q2 30s / q3 : 50s aide voisin / q4 : 1mn 03 aide ⇒3/4 |

## Entretiens et observation

La discussion s'appuie sur les vingt et une phrases soumises aux étudiants. Les commentaires sont construits à la fois à partir des discussions dans les groupes et des observations en classe d'un groupe d'étudiant travaillant sur un des modules de la plateforme. Les textes en italiques sont des transcriptions des discussions des étudiants dans la phase de q-sorting.

## Groupre 1 (vidéo)

| Tout à fait d'accord | Plutôt d'accord | Plutôt pas d'accord | Pas d'accord du tout |
|---|---|---|---|
| **18.** Les QCM finaux sont plus difficiles que les exercices ou le cours<br><br>**20.** C'est impossible de faire les exercices sans écrire sur un papier les calculs. | **2.** Je n'ai jamais regardé l'introduction du module.<br><br>**3.** Les contenus des vidéos sont adaptés à mon niveau en maths.<br><br>**7.** c'est en faisant les exercices que j'ai vraiment compris<br><br>**10.** Le QCM intermédiaire était intéressant pour savoir où j'en étais.<br><br>**17.** Je suis satisfait dès que j'atteins 50% de réussite<br><br>**19.** J'utilise toujours un papier et un crayon pour faire les QCM (finaux et intermédiaires). | **1.** L'introduction de chaque module est intéressante parce que c'est une situation surprenante.<br><br>**5.** il a fallu que je regarde plusieurs fois une (les) vidéos pour comprendre<br><br>**9.** pour les exercices en vidéo, j'ai été obligé de regarder la correction à chaque fois<br><br>**8.** les exercices sont trop difficiles<br><br>**15.** j'aime mieux la correction sur pdf qu'en vidéo<br><br>**21.** Après avoir fait les premiers modules je me suis senti plus armé pour suivre les cours de GC et de GMP | **4.** Les vidéos sont trop longues<br><br>**6.** les exercices corrigés m'ont bien entraîné pour passer les QCM<br><br>**11.** je n'ai fait les QCM intermédiaires que lorsque j'avais raté le QCM final<br><br>**12.** Je préfère commencer par les QCM intermédiaires puis aller voir le cours ou faire les exercices<br><br>**13.** Les exercices supplémentaires sont trop difficiles<br><br>**14.** Les exercices supplémentaires sont trop nombreux et je n'en ai fait qu'une partie<br><br>**16.** je commence toujours par le QCM final |

## Groupe 2 (vidéo)

| Tout à fait d'accord | Plutôt d'accord | Plutôt pas d'accord | Pas d'accord du tout |
|---|---|---|---|
| **4.** Les vidéos sont trop longues<br><br>**14.** Les exercices supplémentaires sont trop nombreux et je n'en ai fait qu'une partie | **3.** Les contenus des vidéos sont adaptés à mon niveau en maths.<br><br>**2.** Je n'ai jamais regardé l'introduction du module.<br><br>**12.** Je préfère commencer par les QCM intermédiaires puis aller voir le cours ou faire les exercices<br><br>**15.** J'aime mieux la correction sur pdf que sur la vidéo | **9.** pour les exercices en vidéo, j'ai été obligé de regarder la correction à chaque fois<br><br>**11.** je n'ai fait les QCM intermédiaires que lorsque j'avais raté le QCM final<br><br>**17.** Je suis satisfait dès que j'atteins 50% de réussite<br><br>**18.** Les QCM finaux sont plus difficiles que les exercices ou le cours<br><br>**19.** J'utilise toujours un papier et un crayon pour faire les QCM (finaux et intermédiaires).<br><br>**20.** C'est impossible de faire les exercices sans écrire sur un papier les calculs<br><br>**21.** Après avoir fait les premiers modules je me suis | **5.** il a fallu que je regarde plusieurs fois les vidéos pour comprendre<br><br>**6.** les exercices corrigés m'ont bien entraîné pour passer les QCM<br><br>**7.** c'est en faisant les exercices que j'ai vraiment compris<br><br>**8.** les exercices sont trop difficiles<br><br>**10.** Le QCM intermédiaire était intéressant pour savoir où j'en étais.<br><br>**13.** Les exercices supplémentaires sont trop difficiles<br><br>**16.** je commence toujours par le QCM final |

| Tout à fait d'accord | Plutôt d'accord | Plutôt pas d'accord | Pas d'accord du tout |
|---|---|---|---|
| | | senti plus armé pour suivre les cours de GC et de GMP | |

## Groupe 3

| Tout à fait d'accord | Plutôt d'accord | Plutôt pas d'accord | Pas d'accord du tout |
|---|---|---|---|
| **6.** les exercices corrigés m'ont bien entraîné pour passer les QCM<br><br>**10.** Le QCM intermédiaire était intéressant pour savoir où j'en étais.<br><br>**14.** Les exercices supplémentaires sont trop nombreux et je n'en ai fait qu'une partie | **4.** Les vidéos sont trop longues<br><br>**8.** les exercices sont trop difficiles<br><br>**9.** pour les exercices en vidéo, j'ai été obligé de regarder la correction à chaque fois | **1.** 1ntroduction de chaque module est intéressante parce que c'est une situation surprenante.<br>**3.** Les contenus des vidéos sont adaptés à mon niveau en maths.<br>**7.** c'est en faisant les exercices que j'ai vraiment compris<br>**13.** Les exercices supplémentaires sont trop difficiles<br>**15.** J'aime mieux la correction sur pdf que sur la vidéo<br>**18.** Les QCM finaux sont plus difficiles que les exercices ou le cours<br>**19.** J'utilise toujours un papier et un crayon pour faire les QCM (finaux et intermédiaires).<br>**21.** Après avoir fait les premiers modules je me suis senti plus armé pour suivre les cours de GC et de GMP | **2.** je n'ai jamais regardé l'introduction du module.<br><br>**5.** il a fallu que je regarde plusieurs fois les vidéos pour comprendre<br><br>**11.** je n'ai fait les QCM intermédiaires que lorsque j'avais raté le QCM final<br><br>**12.** Je préfère commencer par les QCM intermédiaires puis aller voir le cours ou faire les exercices<br><br>**16.** je commence toujours par le QCM final<br><br>**17.** Je suis satisfait dès que j'atteins 50% de réussite<br><br>**20.** C'est impossible de faire les exercices sans écrire sur un papier les calculs. |

## Groupe 4 :

| Tout à fait d'accord | Plutôt d'accord | Plutôt pas d'accord | Pas d'accord du tout |
|---|---|---|---|
| **4.** Les vidéos sont trop longues<br><br>**10.** Le QCM intermédiaire était intéressant pour savoir où j'en étais.<br><br>**14.** Les exercices supplémentaires sont trop nombreux et je n'en ai fait qu'une partie | **3.** Les contenus des vidéos sont adaptés à mon niveau en maths.<br>**6.** les exercices corrigés m'ont bien entraîné pour passer les QCM<br>**7.** c'est en faisant les exercices que j'ai vraiment compris<br>**15.** J'aime mieux la correction sur pdf que sur la vidéo<br>**19.** J'utilise toujours un papier et un crayon pour faire les QCM (finaux et intermédiaires).<br>**20.** C'est impossible de faire les exercices sans écrire sur un papier les calculs.<br>**21.** Après avoir fait les premiers modules je me suis senti plus armé pour suivre les cours de GC et de GMP<br>**12.** Je préfère commencer par les QCM intermédiaires puis aller voir le cours ou faire les | **1.** 1ntroduction de chaque module est intéressante parce que c'est une situation surprenante.<br><br>**5.** il a fallu que je regarde plusieurs fois les vidéos pour comprendre<br><br>**8.** les exercices sont trop difficiles<br><br>**13.** Les exercices supplémentaires sont trop difficiles<br><br>**18.** Les QCM finaux sont plus difficiles que les exercices ou le cours | **2.** je n'ai jamais regardé l'introduction du module.<br><br>**9.** pour les exercices en vidéo, j'ai été obligé de regarder la correction à chaque fois<br><br>**11.** je n'ai fait les QCM intermédiaires que lorsque j'avais raté le QCM final<br><br>**16.** je commence toujours par le QCM final<br><br>**17.** Je suis satisfait dès que j'atteins 50% de réussite |

| | exercices | | |
|---|---|---|---|

Les faits marquants concernant les entretiens et l'activité de q-sorting

**1. L'introduction de chaque module est intéressante parce que c'est une situation surprenante.**

Les groupes semblent unanimes pour ne pas voir l'intérêt de cette introduction, ou même comme certains étudiants l'ont déclaré de ne pas l'avoir vu.
- *Est-ce que ça t'as surpris ?*
- *Ben non, parce que je savais déjà*

**2. je n'ai jamais regardé l'introduction du module.**

Cette deuxième affirmation qui complète et précise la première reçoit des avis divergents et entre en contradiction avec les observations à distance ou en présence. En effet, un groupe se déclare pas d'accord du tout, les deux autres plutôt d'accord.
Les avis sont plutôt contre :
pas « surprenant » : cette vidéo ne les a pas marqués et ne savent pas s'ils l'ont regardée
certains la zappent :  6 élèves passent moins de 20s dessus
pourtant certains restent longtemps dessus : 5 élèves s'attardent plus de 6 mn dessus
KA reste 9 min ( idem pour E) : des arrêts de 30s environ car écrit et fait sur papier en même temps (mais ne le fait pas savoir lors du questionnaire)

*On peut penser que les vidéos de cours et les exercices sont plus attractifs parce que dans la continuité des cours du lycée alors que l'introduction, peut-être pas assez mise en valeur n'apparaît pas nécessaire. Il n'empêche que nous maintenons son importance didactique et qu'une recommandation serait peut-être de la mettre plus en valeur.*

**3. Les contenus des vidéos sont adaptés à mon niveau en maths.**

Dans cette question, deux groupes sont plutôt d'accord et un groupe plutôt pas d'accord.

- *Non trop simple !*

- *T'as fait tous les questionnaires, toi ?*

- *Non, mais les exos sont chauds, mais les vidéos, trop simples.*

- *Oui, mais y'a des trucs, quand même.*

Pour cet étudiant (qui a finalement fait placer cette affirmation plutôt pas d'accord) les exercices sont plus difficiles que les vidéos.

**4. Les vidéos sont trop longues**

Là encore, désaccord entre les groupes. Deux disent être d'accord ou tout à fait d'accord et un groupe pas d'accord du tout.

Il est clair que le ressentiment de la longueur des vidéos dépend du degré de compréhension de la notion qui est présentée. L'observation montre par ailleurs que les étudiants accélèrent la diffusion des vidéos lorsque les contenus leurs paraissent trop simples.

**5. Il a fallu que je regarde plusieurs fois les vidéos pour comprendre**

Les vidéos ne sont *a priori* pas regardées plusieurs fois par les étudiants. Deux interprétations sont possibles : d'une part les vidéos sont suffisamment clairs pour qu'il ne soit pas nécessaire d'y revenir ou d'autre part, les corrigés écrits ou les réponses aux quiz sont plus parlant pour les élèves.

*Malgré le désaccord apparent donné sur les vidéos du cours, il apparaît que les vidéos sont utiles pour les étudiants mais que les utilisations qu'ils en font sont diverses :*

- *lecture continue de la vidéo, ce qui souvent provoque le sentiment de longueur (OC)*

- *lecture avec des arrêts pour vérifier ou faire les exercices (KA)*

- *lecture avec des accélérés au moment des parties supposées connues (PD)*

*Il est donc intéressant de garder ce format de vidéos qui permet à chacun d'utiliser ce support de cours à sa façon.*

**6. Les exercices corrigés m'ont bien entraîné pour passer les QCM**

Deux groupes d'étudiants sont d'accord avec cette affirmation et un pas du tout d'accord. Il est intéressant de creuser un peu plus pour comprendre les raisons de cette dissymétrie. Dans l'entretien de q-sorting, les trois membres de ce groupe déclarent ne pas avoir regardé les exercices avant de faire les QCM. Il s'agit d'une stratégie parfois observée qui consiste dans un premier temps à faire les QCM et de ne revenir aux exercices qu'en cas d'échec. C'est en tout cas ce qui a justifié le classement de ce groupe. Pour d'autres, au contraire, comme on le verra dans les dialogues suivants, le module doit être regardé dans l'ordre proposé. Il ressort donc de l'observation et de l'entretien que les étudiants utilisent les exercices pour renforcer leur compréhension du sujet et se servent des modules de façons différentes. Leur utilité est avérée en ce sens qu'il y a adéquation entre les objectifs des auteurs de la plateforme et des attentes des étudiants.

**7. C'est en faisant les exercices que j'ai vraiment compris le cours**

On retrouve, bien sûr comme précédemment cette idée de la non nécessité des exercices quand le sujet est déjà maîtrisé et de l'intérêt des exercices pour s'entraîner dans le cas contraire, même si la réussite des exercices est mis en balance avec les vidéos

- *Les vidéos, elles ont bien aidés, quand même* déclare un étudiant dans l'entretien.

**8. Les exercices sont trop difficiles**

Le lien est fait par les étudiants entre les exercices et le cours. Il est intéressant de, encore une fois, se rendre compte que la difficulté d'un exercice est intimement lié à la maîtrise du cours comme l'exprime très bien cet étudiant.

- *Moi je mettrais au milieu… Parce que moi si j'ai fait des fautes, c'est parce que j'ai pas compris le cours… C'est sûr, même !*

- *Alors, on le met où*

- *Ben mets le plutôt pas d'accord.*

La possibilité donnée aux étudiants de s'entraîner (exercices) de s'auto-évaluer (Quiz) et de revenir sur la compréhension des concepts ou des techniques en jeu (cours) montre bien l'intérêt de la plateforme dans sa flexibilité d'utilisation. Chaque étudiant, individuellement utilise ou pas les outils mis à sa disposition.

**9. Pour les exercices en vidéo, j'ai été obligé de regarder la correction à chaque fois**

Les réponses montrent que ce n'est pas systématiquement, mais que la correction est bien utile lorsqu'ils jugent le problème trop difficile.

- *Moi, j'étais pas obligé mais je l'ai fait.*

- *Est-ce t'a regardé le corrigé à chaque fois ?*

- *Ben non, pas à chaque fois, mais quand même, quand j'étais pas bien sûr…*

*Dans l'ensemble des réponses, on voit bien l'importance des exercices dans l'entraînement des étudiants. Même lorsque le sujet paraît simple, c'est à travers la réalisation effective de l'exercice qu'ils confortent leur assurance à savoir résoudre des problèmes.*

**10. Le QCM intermédiaire était intéressant pour savoir où j'en étais.**

Un des groupes a classé dans la catégorie "pas du tout d'accord". Cet avis est justifié par le fait que les étudiants de ce groupe allaient directement au quiz final. Pour les autres, cette autoévaluation était tout à fait intéressante pour leur permettre de savoir où ils en étaient dans leurs apprentissages.

- *Oui, moi oui, toujours*, s'écrie une étudiante dès la lecture de l'affirmation.

**11. Je n'ai fait les QCM intermédiaires que lorsque j'avais raté le QCM final**

- *Ah moi je l'ai fait, ça !*

Pourtant, ce groupe d'étudiants classent dans plutôt pas d'accord parce que ce n'est pas à leurs yeux une stratégie générale. D'une façon générale, et ce que l'on confirme dans les questions suivantes (et qui est aussi globalement confirmé par l'observation en classe), les étudiants suivent l'ordre proposé, même si, bien sûr, ils dérogent parfois à la règle, mais, toujours avec une intention d'apprentissage affirmée.

**12. Je préfère commencer par les QCM intermédiaires puis aller voir le cours ou faire les exercices**

Les étudiants s'appuient *a priori* sur la progression proposée par le professeur. Mais on retrouve la flexibilité de la ressource qui permet d'adapter au niveau de connaissance de chaque étudiant l'utilisation qu'il en fait. Le même groupe d'étudiants que celui cité au paragraphe précédent classe cette fois l'affirmation dans "plutôt d'accord" parce que OC l'a déjà fait et l'annonce. Alors que sur cette même affirmation un autre groupe déclare :

- *non, d'abord les exercices*

- *ben oui, dans l'ordre qu'ils t'ont proposé eux*

- *chuis choqué, ils croient qu'on fait pas dans l'ordre*

- *mais peut-être y'en a qui vont faire les exercices et si ils comprennent pas ils vont voir le cours. Pour aller plus vite.*

- *Aaah ! Mais moi, je suis bête et discipliné, on suit l'ordre et voilà.*

- *Ben oui, mais c'est ça…*

Le "y'en a qui" ressemble bien à une situation vécue qui s'oppose à la stratégie prônée par les autres étudiants. Enfin le troisième groupe classe dans pas du tout d'accord :

- *Ben non. A la limite un QCM plus gros, pourquoi pas...*

Ces dialogues montrent encore une fois les différentes appréhensions de la plateforme et des différentes parties des modules.

*Les QCM apparaissent comme un moyen d'auto-évaluation et d'évaluation formative utilisable et utile pour les étudiants. les QCM finaux sont considérés par les étudiants comme des évaluations de la maîtrise des notions rencontrées.*

### 13. Les exercices supplémentaires sont trop difficiles

Assez unanimement, les étudiants pensent qu'ils ne sont pas trop difficiles et bien dans l'esprit du cours complet.

### 14. Les exercices supplémentaires sont trop nombreux et je n'en ai fait qu'une partie

Les trois groupes sont tout à fait d'accord pour dire qu'il y a trop d'exercices supplémentaires et qu'ils n'en ont fait qu'une partie.

*En ce qui concerne les exercices, la profusion dans les modules donne à chaque étudiant la possibilité de s'entraîner en se situant vis-à-vis d'un niveau donné. Même si le grand nombre d'exercices peut sembler superflu, il constitue une base de référence qu'il serait sans doute important de promouvoir auprès des étudiants tout au long de l'année. L'analyse a priori conduite pour construire ces exercices a été satisfaisante en ce sens que les exercices apparaissent acceptables aux étudiants interrogés.*

### 15. J'aime mieux la correction sur pdf que sur la vidéo

Même si la phrase a été classée dans plutôt d'accord ou plutôt pas d'accord par les groupes, les arguments avancés sont assez semblables : la correction en vidéo est bien utile et regardée initialement mais le fait de pouvoir télécharger la correction en vidéo est apprécié. Il ressort des entretiens que la multiplication des supports est un atout pour la plateforme.

- *C'est pas mal non plus aussi en pdf*

- *Ah ben oui*

Et dans un autre groupe

- *Le pdf c'est bien pratique pour résumer mais la vidéo ça te montre bien ce que t'as pas su faire*

**16. Je commence toujours par le QCM final**

Unanimement les étudiants ne valident pas cette stratégie.

- *Non, on suit la démarche* dit un groupe

**17. Je suis satisfait dès que j'atteins 50% de réussite**

Là encore, même si le classement n'est pas le même dans les groupes, les arguments sont comparables et les étudiants déclarent vouloir atteindre plus que 50% de réussite :

- *50% c'est nul*

- *Faut du 200% ! (Rires)*

La conscience que ces exercices doivent être maîtrisés est unanimes même si la possibilité de faire quelques erreurs demeurent parfois présente :

- *50% non, 90% au moins… Bon allez 80% c'est pas mal…*

- *Moi tant que j'ai pas 100% je le refais*

**18. Les QCM finaux sont plus difficiles que les exercices ou le cours**

Un groupe semble considérer que les QCM sont plus difficiles, mais les explications données ne sont pas audibles sur la vidéo. Cependant, les deux autres groupes signalent que les QCM sont dans la ligne des exercices proposés.

- *Non mais moi, tu vois j'ai eu besoin de faire les erreurs pour comprendre. Maintenant le QCM final… C'est quoi la question ?*

- *Ben est-ce que le QCM final est plus dur ?*

- *Non, plutôt pas d'accord.*

Un regard plus approfondi montre que les QCM sont bien adaptés sans difficultés supplémentaires par rapport aux exercices traités. C'est bien sûr une exigence importante pour que l'évaluation puisse être prise en compte par les étudiants comme un indicateur de leur niveau sur le thème du module.

*Il est intéressant de noter que pour les étudiants, la réussite aux QCM est un critère de compréhension de la notion ou de la technique étudiées. En croisant les réponses aux questions 16-17-18 et aux questions 10-11-12, on peut voir que les QCM sont centraux dans la structuration des modules et qu'ils jouent à la fois un rôle d'auto-évaluation dans une perspective d'évaluation formative et d'évaluation sommative.*

**19. J'utilise toujours un papier et un crayon pour faire les QCM (finaux et intermédiaires).**

Dans l'entretien une étudiante signale qu'elle imprime les cours pdf pour les mettre dans son cahier de cours de maths. Mais là encore les stratégies sont différentes et dépendent du niveau de compréhension du thème du module. Les étudiants sont conscients du fait que ce serait nécessaire :

- *Faudrait utiliser, mais nous on le fait pas…*

- *Moi si, je le mettrais là parce que je le fais de temps en temps.*

La réponse donnée dans les groupes portent plus sur le "toujours" que sur la nécessité, bien perçue, d'utiliser papier crayon pour travailler. Les observations en situation montrent bien par ailleurs l'utilisation fréquente de papier et crayon.

Une incitation plus forte dans les modules pourraient sans doute permettre aux étudiants de s'approprier les contenus plus finement.

**20. C'est impossible de faire les exercices sans écrire sur un papier les calculs.**

Les réponses, très contrastées reposent souvent sur des arguments identiques ; dans un groupe qui classe dans "pas du tout d'accord" :

- *Ce n'est pas impossible… Y'a certains calculs, tu peux les faire de tête*

Un groupe lui classe "tout à fait d'accord" signalant que les calculs, dès qu'ils sont un peu difficile, doivent se faire sur papier.

*Ces deux questions et les observations montrent que même si papier-crayon sont utilisés par les étudiants, ils ne constituent pas une forme de capitalisation des travaux réalisés. Il serait sans doute intéressant de façon à ce que le travail ne se perde pas que les étudiants puisse faire une sorte de port-folio du travail réalisé dans les modules en écrivant à la fois*

*leurs succès mais aussi leurs difficultés et la façon de les surmonter. La remarque d'une étudiante signalant qu'elle imprimait les cours pour les garder comme référence dans son cours de maths montre une utilisation qui pourrait être généralisée.*

**21. Après avoir fait les premiers modules je me suis senti plus armé pour suivre les cours de GC et de GMP**

Des avis peu tranchés sur cette question, parfois "plutôt d'accord", parfois "plutôt pas d'accord".

- *Les premiers modules, quand tu vois, priorité de calculs, ça on connaît [...]*

- *Oui mais c'est normal on sort de bac pro, des fois…*

- *Pas du tout d'accord*

- *Ben non, plutôt pas d'accord…*

Un autre groupe :

- *Y'avait des trucs c'était bien de les revoirs, y'avait d'autres trucs, c'était plutôt facile*

- *Alors plutôt d'accord*

Il faudrait bien sûr pour pouvoir donner un avis plus approfondi sur cette question recueillir aussi les points de vue des professeurs. Les avis des étudiants restent très déclaratifs et dans les deux dialogues qui précèdent, on voit apparaître le bien fondé de ces modules (Oui, mais comme on vient de bac pro, des fois… ou Y'avait des trucs, c'était bien de les revoir) même si la réponse donnée n'est pas aussi affirmative que ces avis pourraient le laisser entendre. La discussion qui a suivi montre que ces modules servent à deux niveaux : d'une part, parce qu'ils permettent aux étudiants de se conforter dans les connaissances nécessaires pour la suite de leur intégration en leur donnant un point de référence initial. D'autre part, comme ils l'indiquent dans le dialogue qui suit, les modules servent pour revoir quelques notions qui sont fondamentales dans la suite, même si pour certains, les modules initiaux paraissent très faciles.

- *Ça aide*

- *Ça fait un petit rafraichissement*

- *Même si les premiers modules, c'est trop facile.*

- *C'est peut être un peu trop les bases, mais ça fait pas de mal*

- *Ça serait bien que ça aille un peu plus loin plus proche de ce qu'on fait en méca.*

Un peu plus tard

- *Oui, les fractions, les égalités, oui ça sert*

Dans les questions supplémentaires abordées avec les étudiants, le temps de travail sur chacun des modules a permis aux étudiants de raconter leur façon de travailler. Le temps passé pour chaque module n'excède pas une heure et les modules sont étudiés d'un trait sans interruption.

- *Je devais rester une petite heure.*

- *J'en fais un en entier, puis le deuxième je le fais le lendemain.*

# Conclusion

Les remarques de la première partie de ce document sont bien sûr des conclusions de l'étude et nous voudrions plutôt développer ici les perspectives qui peuvent s'avérer fructueuses pour un développement futur en même temps que de rappeler la construction des contenus de la plateforme. Ce travail est né de la volonté des professeurs de mécanique et de génie civil et de la direction de l'ENEPS de fournir aux étudiants entrant des bases mathématiques leur permettant de suivre les cours sans être perturbés et arrêtés par des difficultés d'ordre mathématique. Le constat initial partait de l'observation empirique d'étudiants qui comprenaient parfaitement les concepts physiques enseignés mais qui ne réussissaient pas du fait de lacunes dans les techniques mathématiques élémentaires. Les difficultés repérées par les enseignants, puis analysées et didactisées ont été à l'origine des différents modules. Les étudiants en thèse de l'IUT de Grenoble ont ainsi en collaboration avec l'équipe de l'IFÉ pensés, réalisés et mis en oeuvre les contenus de ces modules qui sont actuellement disponibles pour les étudiants de l'ENEPS.

L'étude qui débouche sur ce rapport s'est posé initialement les questions de l'utilité, de l'utilisabilité et de l'acceptabilité de cet ensemble de cours pour les étudiants de l'ENEPS. Les observations tant du point de vue ergonomique que didactique ont clairement montré les apports de ces modules de cours pour les étudiants même si quelques éléments sont encore à améliorer comme il est déjà précisé dans la partie "Résultats des observations en 6 remarques". L'expérience de conception et de réalisation peut alors être reproduite pour des modules futurs qui pourraient être développés sur le modèle et avec la méthodologie utilisée dans la conception des six premiers modules. Il est important de garder en tête la nécessité d'un travail approfondie de didactique des mathématiques de façon à offrir aux étudiants un contenu adapté à leurs besoins, scientifiquement irréprochables et utiles pour les enseignements de l'ENEPS.

Les observations réalisées à distance et en présence montrent que les étudiants sont globalement satisfaits de ce qui est proposé dans ces modules et pour la plupart indiquent que les contenus sont des points de repère importants pour pouvoir se situer par rapport aux connaissances attendues en mathématiques pour commencer l'année avec confiance. Le fait de pouvoir se situer, s'auto-évaluer avant les cours est certainement un atout important de la plateforme. Dans les perspectives de travail et d'amélioration, on pourrait imaginer une évaluation diagnostique préalable à l'utilisation des modules. Cette évaluation pourrait dispenser certains étudiants qui maîtrisent les contenus de l'un ou de l'autre des modules et les diriger vers les modules dans lesquels l'évaluation montrerait les difficultés rencontrées. Cette évaluation pourrait certainement encore mieux participer à l'acceptabilité de la plateforme par les étudiants.

# Bibliographie